\documentclass[12pt]{amsart}
\usepackage{amsmath,latexsym,amssymb}

\input xypic
\newtheorem{Theorem}{Theorem}[section]

\newtheorem{Lemma}[Theorem]{Lemma}
\newtheorem{lemma}[Theorem]{Lemma}
\newtheorem{Corollary}[Theorem]{Corollary}
\newtheorem{corollary}[Theorem]{Corollary}

\newtheorem{Conjecture}[Theorem]{Conjecture}

\newtheorem{remark}[Theorem]{Remark}

\def\QED{\hfill$\Box$}

\begin{document}

\title[A Generalization of the Strong Castelnuovo Lemma
]{A Generalization of the Strong Castelnuovo Lemma}

%\author[L. Ghezzi]
%{Laura Ghezzi}

%\thanks{AMS 2000 {\em Mathematics Subject Classification}.
%Primary 13A30; Secondary 13B22, 13H10, 13H15.}

\author{Laura Ghezzi}
\address{Department of Mathematics, New York City College of
Technology-CUNY}
\address{300 Jay Street, Brooklyn, NY 11201,  U.S.A.}
\email{lghezzi@citytech.cuny.edu}

\thanks{This work was supported (in part) by a grant from The City University of New York PSC-CUNY Research Award Program-39}

\thanks {The author would like to thank D. Eisenbud, M.E. Rossi and G. Valla for stimulating conversations about the topic of this paper}

%Abstract for Ottawa: We consider a set $X$ of distinct points in the $n$-dimensional projective space over an algebraically closed field $k$. Let $A$ denote the coordinate ring of $X$, and let $a_i(X)=\dim_k [{\rm
%Tor}_i^R(A,k)]_{i+1}$. Green's Strong Castelnuovo Lemma (SCL) shows that if the points are in general position, then $a_{n-1}(X)\neq 0$ (that is, there are linear syzygies up to order $n-1$) if and only if the points are on a rational normal curve. Cavaliere, Rossi and Valla conjectured that if the points are not necessarily in general position the possible extension of the SCL should be the following: $a_{n-1}(X)\neq 0$ if and only if either the points are on a rational normal curve or in the union of two linear subspaces whose dimensions add up to $n$. In this work we prove the conjecture.

\begin{abstract}
 We consider a set $X$ of distinct points in the $n$-dimensional projective space over an algebraically closed field $k$. Let $A$ denote the coordinate ring of $X$, and let $a_i(X)=\dim_k [{\rm
Tor}_i^R(A,k)]_{i+1}$. Green's Strong Castelnuovo Lemma (SCL) shows that if the points are in general position, then $a_{n-1}(X)\neq 0$
if and only if the points are on a rational normal curve. Cavaliere, Rossi and Valla conjectured in \cite{CRV} that if the points are not necessarily in general position the possible extension of the SCL should be the following: $a_{n-1}(X)\neq 0$ if and only if either the points are on a rational normal curve or in the union of two linear subspaces whose dimensions add up to $n$. In this work we prove the conjecture.
\end{abstract}

\maketitle

%%%%%%%%%%%%%%%%%%%%%%%%%%%%%%%%%%%%%%%%%%%%%%%%%%%%%%%%%%%%%%%%%%%%%%%%%%%%%%%%%%%%%%%%%%%%%%%%%%%%%%%%%%%%%%%%%%%%%%%%%%%%%%%%%%%%%%%%%%%%%%%%%%%%%%%%%%%%%%%%%%%%%%%

\section{Introduction}\label{intro}

Let $k$ be an algebraically closed field, and let
$X=\{P_1,\dots,P_s\}$ be a set of $s\geq n+1$ distinct points in
$\mathbb{P}^n:=\mathbb{P}_k^n$, not contained in any hyperplane.

Let $I=I(X)$ denote the defining ideal of $X$ in the polynomial ring
$R=k[x_0,\dots, x_n]$, and let $A=R/I$ denote its homogeneous coordinate ring.

The graded $R$-module $A$ has a minimal free
resolution
$$0 \longrightarrow F_n  \longrightarrow \dots \longrightarrow F_1 \longrightarrow
R \longrightarrow A\longrightarrow 0,$$
where $F_i=\bigoplus_{j=1}^{\beta_i}R(-d_{ij})$.

Many authors have been interested in the relation between the numerical invariants
of the resolution and the geometric properties of $X$.

We are mostly interested in the ``linear part'' of the resolution,
that is, the syzygies that are determined by linear forms.  This
study has been initiated by Green \cite{Green} and the main idea
coming from his work is that ``a long linear strand in the
resolution has a uniform and simple motivation''. See for example \cite{CRV},
 \cite{EK}, \cite{EK2}, \cite{EGHP}, \cite{EH}, \cite{EP}, \cite{Green2}, \cite{GrLaz}, \cite{P}, \cite{Y}
(this is by no means a complete list) and the literature cited
there.

For every $i=1,\dots, n$, let $a_i:=a_i(X)=\dim_k [{\rm
Tor}_i^R(A,k)]_{i+1}$ denote the multiplicity of the shift $i+1$ in
$F_i$.

It is well known that if $a_i=0$ for some $i$, then $a_j=0$ for all $j\geq i$.
Since $a_1=\dim_k(I_2)$, where $I_2$ denotes the homogeneous part of degree 2 of $I$, we are interested in varieties lying on some quadric.
%maybe some history about $a_i\neq 0$ and $a_n$???

We say that $X$ is in general position if $n+1$ points of $X$ are
never on a hyperplane.

\medskip

A well celebrated result of Green, the Strong Castelnuovo Lemma
(SCL for short), shows that for a set of distinct points in $\mathbb{P}^n$
in general position, we have that $a_{n-1}\neq 0$ (that is, there is
a linear strand of length $n-1$ in the resolution) if and only if the points are on a
rational normal curve of $\mathbb{P}^n$ (see \cite[3.c.6]{Green}).

It is natural to ask what happens if the points are not necessarily in general
position. Cavaliere, Rossi, and Valla conjectured in \cite{CRV} that
the possible extension of the SCL should be the following.

\begin{Conjecture}\label{theconjecture}
For a set $X$ of distinct points spanning $\mathbb{P}^n$, one has
$a_{n-1}\neq 0$ if and only if either the points are on a rational
normal curve or on $\mathbb{P}^k\cup \mathbb{P}^r$ for some positive
integers $k$ and $r$ such that $k+r=n$.
\end{Conjecture}

It follows from \cite[1.2]{CRV} that if the points are on a rational
normal curve or on $\mathbb{P}^k\cup \mathbb{P}^r$ with $k+r=n$, then $a_{n-1}\neq 0$. In view of this result and of the SCL, Conjecture~\ref{theconjecture} can be
restated as follows.

\begin{Conjecture}\label{restated}
If $X$ is not in general position and $a_{n-1}\neq 0$, then $X \subset \mathbb{P}^k\cup \mathbb{P}^r$ for some positive integers $k$ and $r$ with $k+r=n$.
\end{Conjecture}

In this work we prove the following theorem, which appears in Section~\ref{mainresult} as Theorem~\ref{main}.

\begin{Theorem}\label{main2} Let $X$ be a set of distinct points spanning $\mathbb{P}^n$. Fix $i=0,\dots,n-2$. Assume that:
\begin{enumerate}

\item There exist $n-i+1$ points of $X$ on a $\mathbb{P}^{n-i-1}$,

\item $n-i$ points of $X$ are never on a $\mathbb{P}^{n-i-2}$,

\item $a_{n-1}\neq 0$.
\end{enumerate}

Then $X \subset \mathbb{P}^k\cup \mathbb{P}^r$ for some
positive integers $k$ and $r$ such that $k+r=n$.
\end{Theorem}

Notice that if the points are not in general position, then (1) is satisfied for $i=0$. Since (2) is satisfied for $i=n-2$, Theorem~\ref{main2} proves Conjecture~\ref{restated}.

Cavaliere, Rossi, and Valla proved Theorem~\ref{main2} for $i=0$ and $i=1$ (cases they were interested in
for other purposes, see \cite[4.2]{CRV}).

\medskip

Following the philosophy of \cite{CRV}, the main idea of this work is to study explicitly the quadrics passing through the points. We show that there are enough quadrics that ``split'' into the product
of two linear forms to guarantee that $X$ is contained in the union of two linear subspaces whose dimensions add up to $n$.

\medskip

Now we briefly describe the content of this paper. In Section~\ref{prelim} we recall very useful tools from \cite{CRV}. The main point is that $a_{n-1}\neq 0$
 implies that there is at least one nonzero quadric of the form
$F_{abc}=\lambda_{abc} x_ax_b+\mu_{abc} x_ax_c +\nu_{abc} x_bx_c$, $0\leq a<b<c\leq n$, passing through the points. We will refer to such quadrics as ``special quadrics''. Remark~\ref{relation} shows that these special quadrics are ``nicely related''.

The bulk of the paper is given by Section~\ref{generaltheorem}. We prove a general result (Theorem~\ref{general}) showing that if we know that certain special quadrics are reducible, then we can explicitly construct more reducible quadrics passing through the points.
%More precisely we are assuming that there exists $x_j\in \{x_0,\dots,x_n\}$ such that $F_{efj}=x_jL_{ef}$ for all $\{e,f\}$ in a subset of $\{x_0,\dots,x_n\}\setminus \{j\}$. Here $L_{ef}$ is a linear form in the variables $x_e$ and $x_f$, and the main point is to consider the possible monomial or binomial structure of such linear forms.

In Section~\ref{mainresult} we start the proof of Theorem~\ref{main2}. The assumptions guarantee that all quadrics $F_{abj}$, with $\{a,b\}\subset
\{0,\dots, n-i-1\}$ and $j\in \{n-i,\dots, n\}$ ``split'',
$F_{abj}=x_jL_{ab}^j$, where $L_{ab}^j$ is a linear form in $x_a$ and
$x_b$.

Let $W_j$ be the vector space generated by the
linear forms $L_{ab}^j$. First we show that if $W_j=0$ for all $j=n-i,\dots,n$, then $X \subset \mathbb{P}^k\cup \mathbb{P}^r$ for some
positive integers $k$ and $r$ such that $k+r=n$. This statement follows easily from Theorem~\ref{general}.

In Section~\ref{step2} we complete the proof of Theorem~\ref{main2} by proving that if $W_j\neq 0$ for some $j$, then $X \subset \mathbb{P}^k\cup \mathbb{P}^r$ for some positive integers $k$ and $r$ such that $k+r=n$. We first prove the statement when $\dim W_j\geq n-i-1$ (Theorem~\ref{bigdim}). When $\dim W_j< n-i-1$ we use Theorem~\ref{general} as a starting point.

%%%%%%%%%%%%%%%%%%%%%%%%%%%%%%%%%%%%%%%%%%%%%%%%%%%%%%%%%%%%%%%%%%%%%%%%%%%%%%%%%%%%%%%%%%%%%%%%%%%%%%%%%%%%%%%%%%%%%%%%%%%%%%%%%%%%%%%%%%%%%%%%%%%%%%%%%%%%%%%%

\section{Preliminaries}\label{prelim}

In this section we introduce the necessary notation and we recall tools that are very useful in the proof of Theorem \ref{main2}.

\medskip

We compute ${\rm
Tor}_i^R(A,k)$ using a resolution of the field $k$ which can be
obtained from the Koszul complex of $x_0,\dots,x_n$. We fix a
$k$-vector space $V$ of dimension $n+1$. Then the Koszul resolution
of $k$ is given by
$$0\rightarrow \wedge^{n+1}V\otimes R(-n-1)\rightarrow \wedge^{n}V\otimes
R(-n)\rightarrow\dots\rightarrow \wedge V\otimes
R(-1)\rightarrow R\rightarrow k\rightarrow 0.$$ Let
$\delta_i:\wedge^{i}V\otimes R(-i)\rightarrow
\wedge^{i-1}V\otimes R(-i+1)$ be the usual Koszul map. We denote
by $K_{n-2}$ the kernel of $\delta_{n-2}$ in degree $n$. A special
case of \cite[1]{CRV1} gives that
$$a_{n-1}=\dim_k\big[\big(\wedge^{n-2}V\otimes I_2\big)\cap
K_{n-2}\big],$$ where $I_2$ denotes the homogeneous part of degree 2
of the ideal $I$.

\medskip

Let
$e_0,\dots, e_n$ be a $k$-vector basis of $V$. If $j$ is a
$(n-2)$-tuple $\{0\leq j_1<\dots<j_{n-2}\leq n\}$, let
$\epsilon_j:=e_{j_1}\wedge\dots\wedge e_{j_{n-2}}\in
\wedge^{n-2}V$. The following observations play a crucial role.

\begin{remark}\label{alpha}{\rm (\cite[1.3]{CRV}) We have that every element $\alpha \in
\big(\wedge^{n-2}V\otimes I_2\big)\cap K_{n-2}$ can be
written as $\alpha=\sum_{|j|=n-2}\epsilon_j\otimes F_{C_j},$ where
$C_j:=\{0,\dots,n\}\setminus \{j\}$ and $F_{C_j}\in I_2$ is a square
free quadratic form in the variables $x_l,\  l\in C_j$.}\QED
\end{remark}

Therefore if $a_{n-1}\neq 0$ there is at least one nonzero quadric of the form
$$F_{abc}=\lambda_{abc} x_ax_b+\mu_{abc} x_ax_c +\nu_{abc} x_bx_c,$$ $0\leq a<b<c\leq n$, passing through the points.

\begin{remark}\label{relation}{\rm (\cite[1.4]{CRV}) For every
$\{a,b,c,d\}$ such that $0\leq a<b<c<d\leq n$ we
have that
$$(-1)^ax_{a}F_{bcd}+(-1)^{b-1} x_{b}F_{acd}+(-1)^{c-2} x_{c}F_{abd}
+(-1)^{d-3} x_{d}F_{abc}=0. $$}\QED
\end{remark}

%%%%%%%%%%%%%%%%%%%%%%%%%%%%%%%%%%%%%%%%%%%%%%%%%%%%%%%%%%%%%%%%%%%%%%%%%%%%%%%%%%%%%%%%%%%%%%%%%%%%%%%%%%%%%%%%%%%%%%%%%%%%%%%%%%%%%%%%%%%%%%%%%%%%%%%%%%%%%%%%

\section{Reducible quadrics through the points}\label{generaltheorem}

In this section we prove Theorem~\ref{general}. The proof gives an explicit description of certain quadrics passing through the points that ``split'' into the product of two linear forms. Most of the proof of Theorem~\ref{main2} will follow from Theorem~\ref{general}.
% the theorem shows that if certain quadrics split, then more will split.

\medskip

Consider the quadrics $F_{C_j}$ as in Section~\ref{prelim}.
%$$F_{abc}=\alpha_{abc} x_ax_b+\beta_{abc} x_ax_c +\gamma_{abc} x_bx_c,$$
%for $0\leq a<b<c\leq n$.

\begin{Theorem}\label{general}
Let $j\in \{0,\dots,n\}$, and let $\{i_1,\dots,i_m\}\subset\{0,\dots,n\}\setminus\{j\}$ with $i_1<\dots<i_m$.
Suppose that $F_{efj}\in (x_j)$ for all $\{e,f | e\neq f\}\subset\{i_1,\dots,i_m\}$.
Write $F_{efj}=x_jL_{ef}$, and let $V$ be the vector space spanned by the linear forms $L_{ef}$. Let $d:=\dim V$ and suppose that $0<d<m-1$.
Then there exist $t$, with $0<t\leq d$, linear forms $L_1,\dots,L_t$, which are part of a basis of $V$, and linearly independent linear forms $h_1,\dots,h_{m-1-t}$, such that
%in variables not involved in $\{L_1,\dots,L_r\}$,
$$(L_1,\dots,L_t)(h_1,\dots,h_{m-1-t})\subset I.$$
\end{Theorem}

\begin{proof}

%Write general statement about linear forms $L_{ab}=\lambda_{ab}y_a+\mu_{ab}y_b$ with certain properties.

%Assume that there exists a variable $x_j$ such that $F_{abj}=x_jL_{ab}$. Assume that the vector space $V$ spanned by these forms is $d<m-1$, where $m$ is the number %of variables. We rename the variables $y_1,\dots,y_m$.
%proof the rest of this section is devoted to the proof

Let $\mathcal L$ be the set of linear forms $\{L_{ef}| e\neq f\}$. For simplicity of notation we rename the $m$ variables involved in $\mathcal L$ as $y_1,\dots,y_m$. If $L_{ef}\in \mathcal L$ and $e<f$, let $$L_{ef}=\lambda_{ef}y_e+\mu_{ef}y_f.$$
%Let $V_V$ be the variables involved in $V$, and let $V_N=\{y_1,\dots,y_m\}\setminus V_V$.

We may assume that either $j<i_1$, or that $j>i_m$. Applying Remark~\ref{relation} to $\{j,e,f,g\}$ (or $\{e,f,g,j\}$) with $1\leq e<f<g\leq m$ we obtain that

\begin{equation}\label{usefuleq}
F_{efg}=(-1)^{e+j}y_eL_{fg}+(-1)^{f+j-1}y_fL_{eg}+(-1)^{g+j}y_gL_{ef}.
\end{equation}

\medskip

The following lemma will be used often.

\begin{lemma}\label{keylemma} Let $L_{ef}\in \mathcal L$, and suppose that the coefficient of $y_f$ in $L_{ef}$ is not zero. Let $\{u,v\}\subset \{1,\dots,m\}\setminus \{e,f\}$, and let $T$ be a linear form in $y_u$ and $y_v$. Assume that $y_eT\in I$,
and that $L_{eu}$ and $L_{ev}$ are monomials in $y_e$. Then $y_fT\in I$.
\end{lemma}

%\begin{lemma}\label{keylemma} Let $T$ be a linear form in $y_t$ and $y_u$ and assume that $y_eT\in I$. Suppose that the coefficient of $y_f$ in $L_{ef}$ is not zero,
%Let $L_{ef}=\lambda_{ef}y_e+\mu_{ef}y_f$ with $\mu_{ef}\neq 0$,
%and assume that $L_{et}$ and $L_{eu}$ are monomials in $y_e$. Then $y_fT\in I$.
%\end{lemma}

\begin{proof} Suppose that $y_fT\notin I$. Since $y_eT\in I$ there exists a point $E$ such that $y_f(E)\neq 0$, $T(E)\neq 0$ and $y_e(E)=0$. Without loss of generality assume that $y_u(E)\neq 0$, and that $e<f<u$. By (\ref{usefuleq}) $F_{efu}=\pm y_eL_{fu}\pm y_fL_{eu}\pm y_u L_{ef}$. Now $F_{efu}(E)=0$ implies that $\mu_{ef}y_f(E)y_u(E)=0$, a contradiction.
\end{proof}

\medskip

Suppose that at least one among the linear forms in $\mathcal L$ is a non zero monomial, say $L_{ab}$ is a monomial in $y_a$. If the coefficient of $y_c$ in $L_{ac}$ is not zero we say that $y_c$ is {\it connected} to $y_a$ (in one step).

\medskip

We construct inductively  a {\it block of monomials} $B^{ab}$ starting with $L_{ab}$ in the following way. At step 1 we add all the monomials connected to $y_a$. At step $i\geq 2$ we add new monomials connected to the monomials introduced in step $i-1$. In other words, $B^{ab}$ consists of all monomials connected to $y_a$ in a finite number of steps. Notice that the set $B^{ab}$ is part of a basis of $V$, and therefore it contains at most $m-2$ monomials.

\medskip

In what follows $B^{ab}$ denotes the block of monomials starting with $L_{ab}=\lambda_{ab}y_a\neq 0$. We say that $y_a$ is a {\it generator} of $B^{ab}$.

\begin{remark}\

{\rm \begin{enumerate}
%\item If $L_{ab}$ and $L_{ac}$ are non zero monomials in $y_a$, then $B^{ab}=B^{ac}$.
\item If $L_{ab}=\lambda_{ab} y_a\neq 0$ and $L_{ac}=\lambda_{ac} y_a\neq 0$, then $B^{ab}=B^{ac}$.
\item If $y_b\in B^{ab}$ and
%one of the linear forms is a monomial in $y_b$, say
$L_{bc}=\lambda_{bc}y_b\neq 0$, then $B^{ab}=B^{bc}$, since $y_a$ is connected to $y_b$.
\end{enumerate}}\QED
\end{remark}

Next we describe some quadrics that factor into the product of two linear forms.
\medskip

For simplicity of notation let $L_{ab}=L_{12}=\lambda_{12}y_1\neq 0$. By (\ref{usefuleq}) we have that for $s=3,\dots,m$,
$$F_{12s}=(-1)^{j+1}y_1L_{2s}+(-1)^{j+1}y_2L_{1s}+(-1)^{s+j}y_sL_{12}.$$

If $\mu_{1s}=0$ (which is the case if $y_s\notin B^{12}$), then $$ F_{12s}=(-1)^jy_1[((-1)^s\lambda_{12}-\mu_{2s})y_s-(\lambda_{1s}+\lambda_{2s})y_2]=(-1)^jy_1f_s,$$

where \begin{equation}\label{deff} f_s=((-1)^s\lambda_{12}-\mu_{2s})y_s-(\lambda_{1s}+\lambda_{2s})y_2.\end{equation}

\medskip

More generally, applying (\ref{usefuleq}) to $1<u<v$ we have that
$$F_{1uv}=(-1)^{j+1}y_1L_{uv}+(-1)^{u+j-1}y_uL_{1v}+(-1)^{v+j}y_vL_{1u}.$$

If $\mu_{1u}=\mu_{1v}=0$ (which is the case if $y_u,y_v\notin B^{12}$), then $$F_{1uv}=(-1)^jy_1[((-1)^{u-1}\lambda_{1v}-\lambda_{uv})y_u+((-1)^v\lambda_{1u}-\mu_{uv})y_v]=(-1)^jy_1G_{uv},$$
where \begin{equation}\label{defG} G_{uv}=((-1)^{u-1}\lambda_{1v}-\lambda_{uv})y_u+((-1)^v\lambda_{1u}-\mu_{uv})y_v.\end{equation}

In particular, $G_{2v}=f_v$.

\begin{remark}\label{relamongG}{\rm Let $1<u<v<w$. By (\ref{usefuleq}) and (\ref{defG}) we obtain that
$$F_{uvw}=(-1)^{w+j-1}y_w[G_{uv}+(-1)^{u}(\lambda_{1v}+(-1)^{w-1}\mu_{vw})y_u+(-1)^{v-1}(\lambda_{1u}+(-1)^{w-1}\mu_{uw})y_v]$$
$$+(-1)^{v+j-1}\lambda_{uw}y_uy_v+(-1)^{u+j}\lambda_{vw}y_uy_v.$$
In particular, if $\lambda_{uw}=\lambda_{vw}=0$ we have that $F_{uvw}=(-1)^{w+j-1}y_wG_{uv}$ if and only if $\mu_{uw}=(-1)^w\lambda_{1u}$ and $\mu_{vw}=(-1)^w\lambda_{1v}$ if and only if the coefficient of $y_w$ in $G_{uw}$ and in $G_{vw}$ is zero.}\QED\end{remark}

%We say that $B^{ab}$ is a {\it good block} if there exists a non zero monomial $L_{ik}=\lambda_{ik}y_i$ such that $B^{ab}=B^{ik}$ and $y_k\notin B^{ik}$.

\begin{remark}\label{badblock} {\rm Assume that $y_2\in B^{12}$. Let $Y_C=\{y_1,\dots,y_m\}\setminus \{B^{12}\}$ and let $C$ be the set of indexes of the variables in $Y_C$. Notice that $Y_C\neq \emptyset$, since $\dim V<m-1$. If $s\in C$, then $\mu_{1s}=\mu_{2s}=0$. Therefore we have that $y_1f_s\in I$, where $f_s=(-1)^s\lambda_{12}y_s-(\lambda_{1s}+\lambda_{2s})y_2$. If $\lambda_{1s}=\lambda_{2s}=0$ for all $s\in C$, then by Lemma \ref{keylemma} \begin{equation}(B^{12})(Y_C)\subset I\label{eqbadblock},\end{equation} and so the conclusion of Theorem~\ref{general} holds.

If $\lambda_{1s}\neq 0$ for some $s\in C$, then $L_{1s}=\lambda_{1s}y_1\neq 0$ and $B^{12}=B^{1s}$. Notice that $y_s\notin B^{1s}$.
%Since $y_s\notin B^{1s}$, $B^{1s}$ is a good block.
If $\lambda_{2s}\neq 0$ for some $s\in C$, then $L_{2s}=\lambda_{2s}y_2\neq 0$ and $B^{12}=B^{2s}$, since $y_2\in B^{12}$. Notice that $y_s\notin B^{2s}$.
%Since $y_s\notin B^{2s}$, $B^{2s}$ is a good block.

Therefore we may assume that $y_2\notin B^{12}$.}\QED
\end{remark}

%We will assume from now on that $y_b\notin B^{ab}$.

Let $s\neq 2$, and suppose that $y_2,y_s\notin B^{12}$. Then by Lemma \ref{keylemma} \begin{equation}\label{goodblock}(B^{12})f_s\subset I,\end{equation} since for every $y_e\in B^{12}$, $L_{2e}$ and $L_{es}$ are monomials in $y_e$.

More generally, if $y_u,y_v\notin B^{12}$ we have that
\begin{equation}\label{goodblock2}(B^{12})G_{uv}\subset I.\end{equation}

%\begin{remark}\label{oneblock} {\rm Suppose $\{B^{12}\}$ is a basis of $V$. Let $V_N=\{y_1,\dots,y_m\}\setminus \{B^{12}\}\setminus\{y_2\}$ and let $N$ be the set of indexes of the variables in $V_N$. Notice that $V_N\neq \emptyset$, since $\dim V<m-1$. If $s\in N$, then $\mu_{1s}=\mu_{2s}=0$. Therefore we have that $y_1f_s\in I$, where $f_s=(-1)^s\lambda_{12}y_s-(\lambda_{1s}+\lambda_{2s})y_2$. Since $(B^{12})f_s\subset I$, we have that (\ref{conclusion}) holds.}
%\end{remark}

\begin{Lemma}\label{monomial}
Let $L_{ab}=\lambda_{ab}y_a\neq 0$, and assume that the variable $y_b$ does not appear in $V$. Then the conclusion of Theorem~\ref{general} holds.
\end{Lemma}

\begin{proof}
Assume that $L_{ab}=L_{12}=\lambda_{12}y_1$, and construct $B^{12}$. Since $d<m-1$ there exists $y_w\neq y_2$ such that $y_w\notin B^{12}$. Let $Y_C=\{y_1,\dots,y_m\}\setminus\{y_2\}\setminus \{B^{12}\}$, and let $C$ be the set of indexes of the variables in $Y_C$. Let $G_{C}=\{G_{2w}| w\in C\}$, where $G_{2w}=f_w=((-1)^w\lambda_{12}-\mu_{2w})y_w-\lambda_{1w}y_2$. By (\ref{goodblock}) we have that $$(B^{12})(G_{C})\subset I.$$

If for all $w\in C$ the coefficient of $y_w$ in $G_{2w}$ is not zero, then we are done.

Let $M_1=\{w\in C| \mu_{2w}=(-1)^w\lambda_{12}\}$, $Y_{M_1}=\{y_w | w\in M_1\}$ and let $N_1=N\setminus M_1$. We may assume that $M_1\neq\emptyset$. If $w\in M_1$, then $L_{2w}=(-1)^w\lambda_{12}y_w\neq 0$ is a basis element of $V$. Since $d<m-1$ we have that $N_1\neq \emptyset$, and $$(B^{12})(G_{N_1})\subset I,$$ where the set $G_{N_1}=\{G_{2w}| w\in N_1\}$ consists of $|N_1|$ linearly independent linear forms. Up to possibly renaming the variables we may assume that if $w\in M_1$ and $v\in N_1$, then $v<w$.

%$B_1=\{v\in V_1|\  \exists\  w\in M_1\  {\rm with}\  \lambda_{vw}\neq 0\}$, and let $Y_{B_1}=\{y_v | v\in B_1\}$.
Let $Y_{C_1}$ be the set of monomials connected to $Y_{M_1}$ in a finite number of steps, and let $C_1\subset N_1$ be the set of indexes of the variables in $Y_{C_1}$. Then $Y_{C_1}$ is part of a basis of $V$. Let $A_1=N_1\setminus C_1$. Since $d<m-1$ we have that $A_1\neq \emptyset$.
By construction, for all $w\in M_1$ and for all $v\in A_1$ we have that $\lambda_{vw}=0$.
% We modify the original proof (copied below) to have if $t\in M_1$ and $w\in V_1$, then $t>w$, to be consistent with the case w=2.

 By Remark~\ref{relamongG} with $w\in M_1$ and $v\in A_1$, we obtain that $F_{2vw}=(-1)^{w+j-1}y_wG_{2v}$ if and only if the coefficient of $y_w$ in $G_{vw}$ is zero.

\medskip

Let $M_2=\{w\in M_1 | F_{2vw}=(-1)^{w+j-1}y_wG_{2v}\  \forall v\in A_1\}$, $Y_{M_2}=\{y_w | w\in M_2\}$ and let $N_2=M_1\setminus M_2$. If $N_2=\emptyset$, then $(Y_{M_1})(G_{A_1})\subset I$, where $G_{A_1}=\{G_{2v}| v\in A_1\}$. Then by Lemma~\ref{keylemma} we have that $$(B^{12},Y_{M_1}, Y_{C_1})(G_{A_1})\subset I,$$ and the conclusion follows.

So we may assume that $N_2\neq\emptyset$. Let $w\in N_2$. We have that the coefficient of $y_w$ in $G_{vw}$ is not zero, for some $v\in A_1$. Let $G_{N_2}$ be the set of such linear forms $\{G_{vw}\}$. By (\ref{goodblock2}) we have that
$$(B^{12})(G_{N_1},G_{N_2})\subset I,$$
where the set $\{G_{N_1},G_{N_2}\}$ consists of $|N_1|+|N_2|$ linearly independent linear forms.
If $M_2=\emptyset$, then we are done. Otherwise we repeat the procedure.
Let $Y_{C_2}$ be the set of monomials connected to $Y_{M_2}$ in a finite number of steps, and let $C_2\subset N_1\cup N_2$ be the set of indexes of the variables in $Y_{C_2}$.
%$B_2=\{v\in V_2|\  \exists\  w\in M_2\  {\rm with}\  \lambda_{vw}\neq 0\}$, and $Y_{B_2}=\{y_v| v\in B_2\}$.
Then $\{B^{12},Y_{M_2},Y_{C_2}\}$ is part of a basis of $V$. Let $A_2=(N_1\cup N_2)\setminus C_2$ and let $G_{A_2}\subset \{G_{N_1},G_{N_2}\}$ be the set of corresponding linear forms.  Since $d<m-1$ we have that $A_2\neq \emptyset$.
%For all $t\in M_2$ and $w\in A_2$ we have that $\mu_{tw}=0$.

Applying Remark~\ref{relamongG} with $w\in M_2$ and $u,v\in \{A_2\}\cup\{2\}$, we obtain that $F_{uvw}=(-1)^{w+j-1}y_wG_{uv}$ if and only if the coefficient of $y_w$ in $G_{uw}$ and $G_{vw}$ is zero.

Let $M_3=\{w\in M_2 | F_{uvw}=(-1)^{w+j-1}y_wG_{uv}\  \forall G_{uv}\in G_{A_2}\}$, and let $N_3=M_2\setminus M_3$. If $N_3=\emptyset$, then $(Y_{M_2})(G_{A_2})\subset I$. By Lemma~\ref{keylemma} we have that $$(B^{12},Y_{M_2}, Y_{C_2})(G_{A_2})\subset I,$$ and the conclusion follows.

So we may assume that $N_3\neq\emptyset$. Let $w\in N_3$. We have that the coefficient of $y_w$ in $G_{uw}$ or in $G_{vw}$ is not zero, for some $u$ or $v$ in $A_2$. Let $G_{N_3}$ be the set of such linear forms. We have that $$(B^{12})(G_{N_1},G_{N_2},G_{N_3})\subset I.$$
Repeating the argument we obtain that the conclusion of Theorem~\ref{general} is given by \begin{equation}\label{Gak}(B^{12}, Y_{M_k}, Y_{C_k})(G_{A_k})\subset I,\end{equation} for some $k\geq 1$, or by \begin{equation}\label{BGN}(B^{12})(\mathcal G)\subset I,\end{equation} where $\mathcal G$ consists of linear forms $G_{uv}$.
\end{proof}

\begin{Corollary}\label{varnotused} Let $L_{ab}=\lambda_{ab}y_a\neq 0$, let $u\neq b$ and assume that the variable $y_u$ does not appear in $V$. Then either the conclusion of Theorem~\ref{general} holds, or $(B^{ab})y_u\subset I$.
\end{Corollary}

\begin{proof} Let $L_{ab}=L_{12}=\lambda_{12}y_1$. By Remark~\ref{badblock} we may assume that $y_2\notin B^{12}$. Let $u\neq 2$. Then by (\ref{goodblock}) we have that $(B^{12})f_u\subset I$, where $f_u=(-1)^u\lambda_{12}y_u-(\lambda_{1u}+\lambda_{2u})y_2$. If $\lambda_{1u}=\lambda_{2u}=0$, then $(B^{12})y_u\subset I$. If $\lambda_{1u}\neq 0$, then $L_{1u}=\lambda_{1u}y_1\neq 0$, and the conclusion follows from Lemma~\ref{monomial}. If $\lambda_{2u}\neq 0$ we conclude similarly.
%since $y_u$ does not appear in $V$.
%Similarly if $\lambda_{2u}\neq 0$ then $L_{2u}=\lambda_{2u}y_2\neq 0$, and the conclusion follows again from Lemma~\ref{monomial}.
\end{proof}

\medskip

Now we construct part of a basis of $V$ consisting of monomials in the following way. If $L_{a_1b_1}=\lambda_{a_1b_1}y_{a_1}\neq 0$ we construct $B_1=B^{a_1b_1}$ and we assume that $B_1$ is {\it maximal}; that is, it is not contained in any bigger block starting with a monomial in $\mathcal L$.
%In particular $y_{a_1}$ is not connected to any monomial in $\mathcal L$ which is not in $B_1$.
By Remark~\ref{badblock} we also assume that $y_{b_1}\notin B_1$. If among the linear forms in $\mathcal L$ there is a nonzero monomial in one of the remaining variables, say $L_{a_2b_2}=\lambda_{a_2b_2}y_{a_2}\neq 0$ with $y_{a_2}\notin B_1$, we construct $B^{a_2b_2}$ and we assume that it is maximal. We also assume that $y_{b_2}\notin B^{a_2b_2}$. Let $B_2=B^{a_2b_2}\setminus (B^{a_2b_2}\cap B_1)$. Proceeding in this way we construct $B^{a_1b_1},\dots, B^{a_kb_k}$ and $B_1,\dots,B_k$, where $\{B_1,\dots,B_k\}$ is part of a basis of $V$. We have that $y_{b_1}\notin \{B_2,\dots,B_k\}$, otherwise $y_{a_1}$ would be connected to $y_{b_1}$ and $B_1$ would not be maximal.
More generally $y_{b_1},\dots,y_{b_k}\notin \{B_1,\dots, B_k\}$.

\begin{corollary}\label{onlymonomial} If $\{B_1,\dots,B_k\}$, $k\geq 1$, is a basis of $V$, then the conclusion of Theorem~\ref{general} holds.
\end{corollary}

\begin{proof} By construction we have that $y_{b_1}\notin \{B_1,B_2,\dots,B_k\}$, and the conclusion follows from
Lemma~\ref{monomial}.
\end{proof}

The following facts about blocks of monomials will be useful later.

\begin{remark}\label{2blocks}{\rm Let $k\geq 2$, let $B_r,B_s\in \{B_1,\dots,B_k\}$ and assume that $r<s$; that is, $B_s$ has been constructed after $B_r$.
Let $L_{a_ra_s}=\lambda_{a_ra_s}y_{a_r}+\mu_{a_ra_s}y_{a_s}$. We have that $\mu_{a_ra_s}= 0$, otherwise $y_{a_s}\in B_r$. We also have that $\lambda_{a_ra_s}=0$, otherwise $B^{a_rb_r}$ would not be maximal. Therefore $L_{a_ra_s}=0$.
Similarly $L_{a_rw}=0$ for all $y_w\in B_s$.} \QED \end{remark}
%Furthermore, if $y_t\in B_r$, then $L_{tu}$ is a monomial in $y_t$ for all $y_u\in B_s$, otherwise $y_u\in B_r$.

\begin{Corollary}\label{prodblocks} In the notation of Remark~\ref{2blocks}, suppose that $y_{a_r}y_{a_s}\in I$. Then $(B_r)(B_s)\subset I$.
\end{Corollary}
\begin{proof} By Lemma~\ref{keylemma} we have that $(B_s)y_{a_r}\subset I$, since $L_{a_rw}=0$ for all $y_w\in B_s$. Let $y_e\in B_r$. Then $L_{ew}$ is a monomial in $y_e$ for all $y_w\in B_s$, and so applying again Lemma~\ref{keylemma} the conclusion follows.
\end{proof}

\medskip

Now we construct a basis $\{L_1,\dots,L_d\}$ of $V$ consisting of $\{B_1,\dots,B_k\}$ as above and a set $L\subset \mathcal L$ of $l$ binomials (necessarily in variables not involved in $\{B_1,\dots,B_k\}$).
By Corollary~\ref{onlymonomial} we may assume that $l\geq 1$. However, we may have that $k=0$; that is, the basis is given by $L$.

We denote by $V_L$ the variables involved in $L$, and by $V_N$ the variables that do not appear in $V$. Let $|V_L|=s$. Then $|V_N|=m-d+l-s$.
%Let $V_V$ be the variables involved in the basis, $V_{B_i}=\{B_i\}$ the variables involved in $B_i$, $1\leq i\leq k$, $V_L$ the variables involved in $L$, and let $V_N=\{y_1,\dots,y_m\}\setminus V_V$ be the variables that do not appear in $V$.
%Let $\{L_1,\dots,L_d\}=\{B_1,\dots,B_k,L\}$ be the basis of $V$. Let $| L|=l$, $| B_1\cup\dots\cup B_k|=d-l$, and let $| V_L|=s$. Then $| V_N|=m-d+l-s$.

\medskip

Recall that if one among $y_{b_1},\dots,y_{b_k}\in V_N$, then the conclusion of Theorem~\ref{general} follows from Lemma~\ref{monomial}. So we may assume that $y_{b_1},\dots,y_{b_k}\in V_L$, since by construction $y_{b_1},\dots,y_{b_k}\notin \{B_1,\dots, B_k\}$. Then if $y_u\in V_N$, $y_u\neq y_{b_1},\dots,y_{b_k}$, and so by Corollary~\ref{varnotused} we may assume that $(B_1,\dots,B_k)(V_N)\subset I$.

Furthermore, if $L_{ab}\in L$ and $y_u\in V_N$, then by (\ref{usefuleq}) we have $F_{abu}=\pm y_uL_{ab}$, so that $(L)(V_N)\subset I.$
Therefore
\begin{equation}\label{basis+varout}(B_1,\dots,B_k,L)(V_N)\subset I.\end{equation}

If $l\geq s-1$, then the conclusion of Theorem~\ref{general} holds, since $d+(m-d+l-s)\geq m-1$. We may assume that $l\leq s-2$. In particular $l\geq 2$.

%We may also assume that $\{y_{b_1},\dots,y_{b_k}\}\in V_L$, since $\{y_{b_1},\dots,y_{b_k}\}\in V_V$ and $y_{b_1},\dots,y_{b_k}\notin B_1\cup\dots\cup B_k$, otherwise the blocks would not be maximal.

\medskip

Next we need some facts about the binomials in $L$. If $C\subset \mathcal L$, we denote by $V_C$ the set of variables involved in $C$.

\begin{lemma}\label{blocksofbinomials} Suppose that $\{L_1,\dots,L_u\}\subset \mathcal L$ are linearly independent binomials in $q$ variables. Then $\{L_1,\dots,L_u\}$ is a disjoint union of subsets, $\{L_1,\dots,L_u\}=A\cup D_1\cup\dots\cup D_p$, with $p=q-u$, such that $|V_{D_j}|=|D_j| +1$ for all $j=1,\dots,p$, and $|V_{A}|=|A|$.
\end{lemma}
\begin{proof}
%Group binomials together in the following way.
We start with $L_{ab}\in \{L_1,\dots,L_u\}$. At step 1 we add the binomials containing the variables $y_a$ or $y_b$. At step $i\geq 2$ we add binomials containing variables introduced in the previous step. Since the binomials are linearly independent, the subset $S$ thus constructed has the property that either $|V_S|=|S| +1$, or $|V_S|=|S|$. Now if $\{L_1,\dots,L_u\}\setminus S\neq \emptyset$ we repeat the procedure.

Let $D_j$, $j=1,\dots,p$, be the subsets with $|V_{D_j}|=|D_j| +1$, and let $A$ be the union of the remaining subsets. We have that $u=|A|+|D_1|+\dots+|D_p|$ and $q=|A|+(|D_1|+1)+\dots+(|D_p|+1)$, so that $q-u=p$.
\end{proof}

By Lemma~\ref{blocksofbinomials} we have that $L=D_0\cup D_1\cup\dots\cup D_p$, where $p=s-l\geq 2$, $|V_{D_j}|=|D_j| +1$ for all $j=1,\dots,p$, and $|V_{D_0}|=|D_0|$. We may have $D_0=\emptyset$.

\begin{lemma}\label{zerobinomials} Let $L_{ef}\in \mathcal L$. In the above set-up, we have that $L_{ef}=0$ if $y_e\in V_{D_j}$ for some $j=1,\dots,p$, and $y_f\in V_L\setminus V_{D_j}$.
\end{lemma}
\begin{proof} For simplicity of notation, let $L=\{L_1,\dots,L_l\}$ and $D_j=\{L_1,\dots,L_r\}$. Write $L_{ef}=\alpha_1L_1+\dots+\alpha_rL_r+\alpha_{r+1}L_{r+1}+\dots+\alpha_lL_l$. It suffices to show that $\alpha_1=\dots=\alpha_r=0$. By construction $L_1=L_{i_1i_2}$ is a binomial in two variables $y_{i_1}$ and $y_{i_2}$, and for $2\leq m\leq r$, $L_m=L_{a_mi_{m+1}}$ is a binomial in $y_{a_m}$ and $y_{i_{m+1}}$, where $a_m\in \{i_1,\dots,i_m\}$. Then $e=i_c$ for some $1\leq c\leq r+1$. It follows that
\begin{equation}\label{bin} L_{ef}=\alpha_1L_{i_1i_2}+\alpha_2L_{a_2i_3}+\dots+\alpha_{c-1}L_{a_{c-1}i_c}+\alpha_{r+1}L_{r+1}+\dots+\alpha_lL_l, \end{equation}
and that $L_{ef}=0$, if $c=1,2$.
Let $c\geq 3$. There exists $i_k\in \{i_1,\dots,i_{c-1}\}$ such that $i_k\neq a_2,\dots, a_{c-1}$. Therefore $y_{i_k}$ appears in (\ref{bin}) only once; in $L_{a_{k-1}i_k}$ if $k\geq 3$, or in $L_{i_1i_2}$ if $k=1,2$. So if $k\geq 2$ we have that $\alpha_{k-1}=0$, and if $k=1$ we have that $\alpha_1=0$. Now we repeat the procedure to obtain that $\alpha_1=\dots=\alpha_r=0$.
\end{proof}

It follows from Lemma~\ref{zerobinomials} that $L_{ac}=0$ if $y_a\in V_{D_j}$ for some $j=0,\dots,p$, and $y_c\in V_L\setminus V_{D_j}$. Hence if $y_a,y_b\in V_{D_j}$ and $y_c\in V_L\setminus V_{D_j}$, by (\ref{usefuleq}) we have that $F_{abc}=\pm y_cL_{ab}$. Therefore for all $j=0,\dots,p$,
\begin{equation}\label{2setsofbin}(D_j)(V_L\setminus V_{D_j})\subset I.\end{equation}

\begin{remark}\label{onlybinomials}{\rm Suppose that the basis of $V$ consists only of binomials. We have that \begin{equation}\label{eqonlybin}(D_1)(V_L\setminus V_{D_1},V_N)\subset I.\end{equation} Let $|D_1|=r$.  The conclusion of Theorem~\ref{general} holds, since $r+(s-r-1)+(m-s)=m-1$.}
%Therefore we may assume that the basis of $V$ contains a monomial.
\end{remark}

\begin{lemma}\label{monandbin} Let $\{B_1,\dots,B_k,L\}$ be a basis of $V$ and suppose that $L=U\cup Z$ where $L_{uz}=0$ if $y_u\in V_U$ and $y_z\in V_Z$. Let $B\in \{B_1,\dots,B_k\}$ be obtained from $L_{ab}=\lambda_{ab}y_{a}\neq 0$ and assume that $y_{b}\in V_U$. Suppose that there exists $y_t\in V_Z$ with $\lambda_{at}\neq 0$. Then the conclusion of Theorem~\ref{general} holds.
\end{lemma}

\begin{proof}
Assume that $L_{ab}=L_{12}=\lambda_{12}y_1\neq 0$. Then $y_2\in V_U$ and there exists $y_t\in V_Z$ such that $\lambda_{1t}\neq 0$.

Let $s$ be such that $y_s\in V_N\cup V_Z$. Then by (\ref{deff}) $f_s=(-1)^s\lambda_{12}y_s-\lambda_{1s}y_2$, since $L_{2s}=0$. Let $f_{N\cup Z}$ be the set of such linear forms.

%Let $s$ be such that $y_s\in V_Z$. Then $L_{2s}=0$, since $y_2\in V_U$ and $y_s\in V_Z$, and so again $(B_1)f_s\subset I$, where $f_s=(-1)^s\lambda_{12}y_s-\lambda_{1s}y_2$.

Let $s$ be such that $y_s\in V_U$ and $y_s\neq y_2$. Then by (\ref{defG}) $G_{st}=(-1)^{s-1}\lambda_{1t}y_s+(-1)^t\lambda_{1s}y_t$, since $L_{st}=0$. Let $G_{U}$ be the set of such linear forms.

Let $V_C=\{B_1,\dots,B_k\}\setminus\{B^{ab}\}$, and let $C$ be the set of indexes of the variables in $V_C$. Let $f_{C}=\{f_{s}| s\in C\}$, where $f_s=((-1)^s\lambda_{12}-\mu_{2s})y_s-\lambda_{1s}y_2$, since $\lambda_{2s}=0$. By (\ref{goodblock2}) we have that $$(B^{12})(f_{C},f_{N\cup Z},G_{U})\subset I,$$
where $\{f_{N\cup Z},G_{U}\}$ are linearly independent.

If for all $s\in C$ the coefficient of $y_s$ in $f_{s}$ is not zero, then we are done.

Otherwise let $M_1=\{s\in C| \mu_{2s}=(-1)^s\lambda_{12}\}$ and we proceed as in the proof of Lemma~\ref{monomial}.
%We have that $$(B_1)(f_{B_2},\dots,f_{B_k},f_{V_N},f_{V_Z},G_{V_U\setminus \{y_2\},t})\subset I,$$ and $f_{V_N},f_{V_Z},G_{V_U\setminus \{y_2\},t}$ are linearly independent.
%Let $y_s\in B_2\cup\dots\cup B_k$. We have $f_s=((-1)^s\lambda_{12}-\mu_{2s})y_s-\lambda_{1s}y_2$, since $\lambda_{2s}=0$. If the coefficient of $y_s$ in $f_s$ is not zero for all $y_s\in B_2\cup\dots\cup B_k$ then the conclusion follows.
\end{proof}

\begin{corollary}\label{mon+bin} Let $B\in \{B_1,\dots,B_k\}$ be obtained from $L_{ab}=\lambda_{ab}y_a\neq 0$ and assume that $y_b\in V_{D_j}$ for some $j=0,\dots,p$. Then either the conclusion of Theorem~\ref{general} holds, or $(B^{ab})(V_L\setminus V_{D_j})\subset I$.
% (and so $(B)(V_L\setminus V_{D_j})\subset I$).
\end{corollary}
\begin{proof} Let $y_c\in V_L\setminus V_{D_j}$. Then by (\ref{usefuleq}) $F_{abc}=\pm y_a L_{bc}\pm y_b L_{ac}\pm y_c L_{ab}=\pm \lambda_{ac}y_ay_b\pm \lambda_{ab}y_ay_c$. If $\lambda_{ac}\neq 0$ the conclusion follows by Lemma~\ref{monandbin}. If $\lambda_{ac}=0$, then $y_ay_c\in I$, and so by Lemma~\ref{keylemma} $(B^{ab})y_c\subset I$.
\end{proof}

\medskip

Back to the proof of Theorem~\ref{general}, recall that a basis of $V$ is given by $$\{B_1,\dots,B_k,D_0,D_1,\dots,D_p\},$$ with $p\geq 2$ and $k\geq 1$, by Remark~\ref{onlybinomials}. We have $L=D_0\cup D_1\cup\dots\cup D_p$, $|L|=l$, $|V_L|=s$. Let $|V_{D_j}|=s_j$, for $0\leq j\leq p$. Then $|D_0|=s_0$ and $|D_j|=s_j-1$ for $1\leq j\leq p$. Let $q\in \{1,\dots,k\}$. If $B_q$ is obtained from $L_{a_qb_q}=\lambda_{a_qb_q}y_{a_q}\neq 0$, we have that $y_{b_q}\in V_{D_{j_q}}$ for some $j_q\in \{0,\dots,p\}$.

\medskip

By (\ref{basis+varout}), (\ref{2setsofbin}) and Corollary~\ref{mon+bin} we may assume that \begin{equation}\label{onemonblock}(B_1,D_{j_1})(V_N,V_L\setminus V_{D_{j_1}})\subset I.\end{equation} If $k=1$ the conclusion of Theorem~\ref{general} holds, so we may assume that $k\geq 2$.

\medskip

We divide $\{B_2,\dots,B_k\}$ in two groups. Let $\{B_{d_2},\dots,B_{d_r}\}$ be such that there exist $y_{t_2},\dots,y_{t_r}\in V_{D_{j_1}}$ such that $L_{a_{d_2}t_2},\dots,L_{a_{d_r}t_r}$ are nonzero monomials in $y_{a_{d_2}},\dots,y_{a_{d_r}}$, respectively.
%there exists a choice of $y_{b_{d_2}},\dots,y_{b_{d_p}}$ with $y_{b_{d_2}},\dots,y_{b_{d_p}}\in V_{D_{j_1}}$.
Let $\{B_{d_{r+1}},\dots,B_{d_k}\}$ be the remaining blocks. Then for every $y_t\in V_{D_{j_1}}$ we have that $L_{a_{d_{r+1}}t}=\dots=L_{a_{d_k}t}=0$. We may assume that
\begin{equation}\label{almostfinal}(B_1,B_{d_2},\dots,B_{d_r},D_{j_1})(V_N,V_L\setminus V_{D_{j_1}})\subset I.\end{equation}

\begin{Corollary}\label{interactingblocks} In the above notation, assume that $\bar B_1\in \{B_1,B_{d_2},\dots,B_{d_r}\}$ and $\bar B_2\in \{B_{d_{r+1}},\dots,B_{d_k}\}$. Then $(\bar B_1)(\bar B_2)\subset I$.
\end{Corollary}
\begin{proof} Assume that $\bar B_1$ is obtained from $L_{cd}=\lambda_{cd}y_c$ with $y_d\in V_{D_{j_1}}$, and $\bar B_2$ is obtained from $L_{ef}=\lambda_{ef}y_e$. By construction we have that  $L_{de}=0$, and by Remark~\ref{2blocks} we have that $L_{ce}=0$. Then by (\ref{usefuleq}), $F_{cde}=\pm \lambda_{cd}y_cy_e$. Therefore $y_cy_e\in I$, and the conclusion follows from Corollary~\ref{prodblocks}.
\end{proof}

\begin{Corollary}\label{BD} In the above notation, assume that $\bar B_2\in \{B_{d_{r+1}},\dots,B_{d_k}\}$. Then $(\bar B_2)(D_{j_1})\subset I$.
\end{Corollary}
\begin{proof} Assume that $\bar B_2$ is obtained from $L_{ef}=\lambda_{ef}y_e$. Let $L_{ab}\in D_{j_1}$. Then $L_{ae}=L_{be}=0$, and so by (\ref{usefuleq}) $F_{abe}=\pm y_eL_{ab}$. Since $y_e(D_{j_1})\subset I$, by Lemma~\ref{keylemma} we have that $(\bar B_2)(D_{j_1})\subset I$.
\end{proof}

\medskip

Finally, by (\ref{almostfinal}), Corollary~\ref{interactingblocks} and Corollary~\ref{BD}, we have that \begin{equation}\label{final}(B_1,B_{d_2},\dots,B_{d_r},D_{j_1})(V_N,V_L\setminus V_{D_{j_1}},B_{d_{r+1}},\dots,B_{d_k})\subset I.\end{equation}

Now, if $j_1\neq 0$ the conclusion of Theorem~\ref{general} follows since $(d-l)+(s_{j_1}-1)+(m-d+l-s)+(s-s_{j_1})=m-1$. Similarly, if $j_1= 0$ the conclusion follows since $(d-l)+s_{0}+(m-d+l-s)+s-s_{0}=m$.

\end{proof}

%%%%%%%%%%%%%%%%%%%%%%%%%%%%%%%%%%%%%%%%%%%%%%%%%%%%%%%%%%%%%%%%%%%%%%%%%%%%%%%%%%%%%%%%%%%%%%%%%%%%%%%%%%%%%%%%%%%%%%%%%%%%%%%%%%%%%%%%%%%%%%%%%%%%%%%%%%%%%%%%%%%%%%

\section{Main Result}\label{mainresult}

In this section we start the proof of the main theorem, stated in Section~\ref{intro} as Theorem~\ref{main2}.

\begin{Theorem}\label{main} Let $X$ be a set of distinct points spanning $\mathbb{P}^n$. Fix $i=0,\dots,n-2$. Assume that:
\begin{enumerate}

\item There exist $n-i+1$ points of $X$ on a $\mathbb{P}^{n-i-1}$,

\item $n-i$ points of $X$ are never on a $\mathbb{P}^{n-i-2}$,

\item $a_{n-1}\neq 0$.
\end{enumerate}

Then $X \subset \mathbb{P}^k\cup \mathbb{P}^r$ for some
positive integers $k$ and $r$ such that $k+r=n$.
\end{Theorem}

%Since (1) is satisfied for $i=0$ and (2) is satisfied for $i=n-2$, the Conjecture follows.

\begin{proof}

Since $n-i+1$ points are in $\mathbb{P}^{n-i-1}$ they must span it, otherwise we get $n-i+1$ points on $\mathbb{P}^{n-i-2}$. After a change of coordinates we may assume that the coordinate points are on $X$ and that the linear space $x_{n-i}=x_{n-i+1}=\dots=x_n=0$ contains $n-i+1$ points of $X$. This linear space contains $n-i$ coordinate points, so it contains an ``extra point'' $Q=(q_0,\dots,q_{n-i-1},0,\dots,0)$. Notice that $q_l\neq 0$ for all $l=0,\dots,n-i-1$, otherwise we would have $n-i$ points in $\mathbb{P}^{n-i-2}$.

\medskip

Let $0\leq a<b\leq n-i-1$ and $n-i\leq j\leq n$. Consider the quadrics
$$F_{abj}=\lambda_{abj} x_ax_b+\mu_{abj} x_ax_j +\nu_{abj} x_bx_j,$$ defined in Section~\ref{prelim}.
Since $F_{abj}(Q)=0$, we have that $$F_{abj}=x_jL^j_{ab},$$ where $L^j_{ab}$ is a linear form in $x_a$ and $x_b$.

Using Remark~\ref{relation} we obtain the following result, which is a generalization of Claim 5 of \cite{CRV}.

\begin{lemma}\label{relquad}\

\begin{enumerate}

\item Let $0\leq a<b<c\leq n-i-1$ and $n-i\leq j\leq n$.
Then $$F_{abc}=(-1)^{a+j}x_aL^j_{bc}+(-1)^{b+j-1}x_bL^j_{ac}+(-1)^{c+j}x_cL^j_{ab}.$$

\item Fix $d$ and $e$ such that $n-i\leq d<e\leq n$. Then there exists $P_{de}=\lambda_{de}x_d+\mu_{de}x_e$ such that for all $s=0,\dots,n-i-1$,
$$F_{sde}=(-1)^sx_sP_{de}+\nu_{sde}x_dx_e.$$
%where $\alpha_{t,\{d,e\}}\in k$.
    %($P_{st}$ does not depend on $h$).

\item Let $0\leq a<b\leq n-i-1$ and $n-i\leq d<e\leq n$. We have that
$$(-1)^ax_a\nu_{bde}+(-1)^{b-1}x_b\nu_{ade}+(-1)^{d-2}L^{e}_{ab}+(-1)^{e-3}L^{d}_{ab}=0.$$

\end{enumerate}
\end{lemma}

\medskip

Let $j\in \{n-i,\dots,n\}$ and let $W_j$ be the $k$-vector space generated by the linear forms $L^j_{ab}$.

\medskip

The proof of Theorem~\ref{main} consists of two main steps:

\begin{enumerate}

\item If $W_j=0$ for all $j=n-i,\dots,n$, then $X \subset \mathbb{P}^k\cup \mathbb{P}^r$ for some
positive integers $k$ and $r$ such that $k+r=n$.

\item If $W_j\neq 0$ for some $j\in \{n-i,\dots,n\}$, then $X \subset \mathbb{P}^k\cup \mathbb{P}^r$ for some
positive integers $k$ and $r$ such that $k+r=n$.

\end{enumerate}

\subsection{Proof of Step 1}

Suppose that $W_j=0$ for all $j=n-i,\dots,n$. Then $F_{abc}=0$, if $0\leq a<b\leq n-i-1$ and $n-i\leq c\leq n$. If $0\leq a<b<c\leq n-i-1$, we also have that $F_{abc}=0$, by Lemma~\ref{relquad} (1). In particular, if $i=0$ we have that $\alpha=0$, a contradiction to Remark~\ref{alpha}, since $a_{n-1}\neq 0$. So we may assume that $i\geq 1$.

If $0\leq a\leq n-i-1$ and $n-i\leq b<c\leq n$, by Lemma~\ref{relquad} (2) and (3) we have that
%$F_{abc}=(-1)^ax_aP_{bc}+\alpha_{a,\{b,c\}}x_bx_c$. Fix $d\in\{0,\dots,n-i-1\}$ with $d\neq a$. By Lemma~\ref{relquad} (3) we have that $\alpha_{a,\{b,c\}}=0$,
$$F_{abc}=(-1)^ax_aP_{bc}.$$

Therefore, if $P_{bc}=0$ for all $n-i\leq b<c\leq n$, we have $F_{abc}=0$ for all $0\leq a\leq n-i-1$ and all $n-i\leq b<c\leq n$. It follows from Remark~\ref{relation} that $F_{abc}=0$ for all $n-i\leq a< b<c\leq n$. Hence $\alpha=0$, a contradiction.

\medskip

Let $W$ be the $k$-vector space generated by the linear forms $P_{bc}$ with $n-i\leq b<c\leq n$, and let $d:=\dim W>0$.

We have that $$(x_0,\dots,x_{n-i-1})W\subset I.$$

We may assume that $d<i$, otherwise $X\subset \mathbb{P}^{i}\cup \mathbb{P}^{n-i}$ and we conclude. Then by Theorem~\ref{general} there exist $t$, with $0<t\leq d$, linear forms $L_1,\dots,L_t$, which are part of a basis of $W$, and linearly independent linear forms $h_1,\dots,h_{i-t}$ such that $$(L_1,\dots,L_t)(h_1,\dots,h_{i-t},x_0,\dots,x_{n-i-1})\subset I.$$
Hence $X\subset \mathbb{P}^{n-t}\cup \mathbb{P}^{t}$ and the conclusion of Theorem~\ref{main} holds. This concludes the proof of Step 1.

\section{Proof of the Main Theorem}\label{step2}
%Here we assume that j=n, so that notation is consistent and we do not have to define f_{tz} and $F_{wz} separately.

In this section we complete the proof of Theorem~\ref{main} by proving Step 2.

\medskip

Fix $j\in \{n-i,\dots,n\}$ such that $W_j\neq 0$. Recall that $W_j$ is the $k$-vector space generated by the linear forms $L^j_{ab}$, where $F_{abj}=x_jL^j_{ab},$ for $0\leq a<b\leq n-i-1$.

%We have that $\dim W_j\leq n-i$.

Let ${\mathcal L}^j$ denote the set of linear forms $\{L^j_{ab}\mid 0\leq a<b\leq n-i-1\}$.

If $0\leq d<e<f\leq n$, we have $$F_{def}=\lambda_{def}x_dx_e+\mu_{def}x_dx_f+\nu_{def}x_ex_f.$$

\begin{remark}\label{defg}{\rm Let $0\leq d<e<f<g\leq n$. By Remark~\ref{relation} we have the following equations:
\begin{equation}\label{defg1} (-1)^d\lambda_{efg}+(-1)^{e-1}\lambda_{dfg}+(-1)^{f-2}\lambda_{deg}=0.\end{equation}
\begin{equation}\label{defg2} (-1)^d\mu_{efg}+(-1)^{e-1}\mu_{dfg}+(-1)^{g-3}\lambda_{def}=0.\end{equation}
\begin{equation}\label{defg3} (-1)^d\nu_{efg}+(-1)^{f-2}\mu_{deg}+(-1)^{g-3}\mu_{def}=0.\end{equation}
\begin{equation}\label{defg4} (-1)^{e-1}\nu_{dfg}+(-1)^{f-2}\nu_{deg}+(-1)^{g-3}\nu_{def}=0.\end{equation}\QED}
\end{remark}

In what follows we give an explicit description of quadrics that are multiple of $x_j$. \medskip

Let $\{d,e | d\neq e \}\subset \{0,\dots,n\}$, and assume that $F_{dej}\in (x_j)$. Without loss of generality, assume that $d<e<j$. Then $\lambda_{dej}=0$, and
$$F_{dej}=x_j(\mu_{dej}x_d+\nu_{dej}x_e)=x_jH_{de},$$ where
$$H_{de}=\mu_{dej}x_d+\nu_{dej}x_e.$$
In particular, if $0\leq d<e\leq n-i-1$, $H_{de}=L^j_{de}$.
%unify notation?

\begin{Lemma}\label{xgGxfG} Let  $\{d,e,f,g\}\subset \{0,\dots,n\}$. Let $T_{de}$ be a linear form in $x_d$ and $x_e$. Assume that $x_gT_{de}\in I$, $x_fT_{de}\notin I$, and that $F_{deg}\in (x_g)$.
Then $F_{dfg}\in (x_g)$ and $F_{efg}\in (x_g)$.
\end{Lemma}
\begin{proof} There exists a point $E$ such that $x_f(E)\neq 0$, $T_{de}(E)\neq 0$, and $x_g(E)=0$. We may assume that $d<e<f<g$.
We have $\lambda_{deg}=0$.
If $x_d(E)\neq 0$, $F_{dfg}(E)=0$
implies that $\lambda_{dfg}=0$.
Similarly if $x_e(E)\neq 0$, then
$\lambda_{efg}=0$. Now (\ref{defg1}) implies that $\lambda_{dfg}=\lambda_{efg}=0$.
\end{proof}

\begin{Corollary}\label{xjGxfG} Let $\{d,e,f\}\subset \{0,\dots,n\}$. Assume that $F_{dej}=x_jH_{de}$ and $x_fH_{de}\notin I$.
Then $F_{dfj}=x_jH_{df}$ and $F_{efj}=x_jH_{ef}$. Furthermore, if $F_{def}\in (x_f)$, we have that the coefficient of $x_f$ in $H_{df}$ or in $H_{ef}$ is not zero.
\end{Corollary}

\begin{proof}
%We have that $F_{dej}=x_j(\mu_{dej}x_d+\nu_{dej}x_e)$.
The first assertion follows from Lemma~\ref{xgGxfG}.
We may assume that $d<e<f<j$. We have that $F_{dfj}=x_j(\mu_{dfj}x_d+\nu_{dfj}x_f)$, $F_{efj}=x_j(\mu_{efj}x_e+\nu_{efj}x_f)$ and $\lambda_{def}=0$.
 If $(\nu_{dfj},\nu_{efj})=(0,0)$,
by (\ref{defg3}) and (\ref{defg4}) we have that $\mu_{def}=(-1)^{f-j}\mu_{dej}$, and $\nu_{def}=(-1)^{f-j}\nu_{dej}$. Then $F_{def}=(-1)^{f-j}x_fH_{de}$,
a contradiction, since $x_fH_{de}\notin I$.
\end{proof}

\begin{Theorem}\label{bigdim} If $\dim W_j\geq n-i-1$, then $X \subset \mathbb{P}^k\cup \mathbb{P}^r$ for some
positive integers $k$ and $r$ such that $k+r=n$.
% so that the conclusion of Theorem~\ref{main} holds.
\end{Theorem}

\begin{proof} We have that $x_jW_j\subset I$. If $(x_{n-i},\dots,x_n)W_j\subset I$, the conclusion follows. Let $A_0=\{x_{n-i},\dots,x_n\}$, and let
$V_1=\{x_u\in A_0 |\ x_uL^j_{ab}\notin I\ {\rm for\  some}\ L^j_{ab}\in W_j \}$. We may assume that $V_1\neq \emptyset$.
Since $F_{abu}=x_uL^u_{ab}$, by Corollary~\ref{xjGxfG} we have that
$x_u\in V_1$ yields a linear form $H_{c_1u}$, $c_1\in \{a,b\}$, with coefficient of $x_u$ different from zero, such that $F_{c_1uj}=x_jH_{c_1u}$.
Let $H_{V_1}$ be the set of such linear forms $H_{c_1u}$.
We have that $x_j(W_j,H_{V_1})\subset I.$ Let $A_1=A_0\setminus V_1$.
If $A_1=\{x_j\}$, then $X \subset \mathbb{P}^{n-1}\cup \mathbb{P}^1$, so we may assume that $\{x_j\}\subsetneq A_1$.
If $(A_1)(H_{V_1})\subset I$, the conclusion follows.

Let $V_2=\{x_v\in A_1 |\ x_vH_{c_1u}\notin I\ {\rm for\  some}\ H_{c_1u}\in H_{V_1} \}$. Since $x_uL^j_{ab}\notin I$, and $x_vL^j_{ab}\in I$, by Lemma~\ref{xgGxfG}
we have that $F_{c_1uv}\in (x_v)$.
By Corollary~\ref{xjGxfG},
$x_v\in V_2$ yields a linear form $H_{c_2v}$, $c_2\in \{c_1,u\}$, with coefficient of $x_v$ different from zero, such that $F_{c_2vj}=x_jH_{c_2v}$.
Let $H_{V_2}$ be the set of such linear forms.

Let $l\geq 2$, and $A_l=A_{l-1}\setminus V_l$. We may assume that
$$x_j(W_j,H_{V_1},\dots, H_{V_l})\subset I,$$ and $$(A_l)(W_j,H_{V_1},\dots, H_{V_{l-1}})\subset I.$$ If $(A_l)(H_{V_l})\subset I$, the conclusion follows.

Let $V_{l+1}=\{x_z\in A_l |\ x_zH_{c_lw}\notin I\ {\rm for\  some}\ H_{c_lw}\in H_{V_l} \}$. By inductively applying Lemma~\ref{xgGxfG}
%maybe state the induction step
we have that $F_{c_lwz}\in (x_z)$.
By Corollary~\ref{xjGxfG},
$x_z$ yields a linear form $H_{c_{l+1}z}$, $c_{l+1}\in \{c_l,w\}$, with coefficient of $x_z$ different from zero, such that $F_{c_{l+1}zj}=x_jH_{c_{l+1}z}$.
Repeating the procedure we obtain the desired conclusion.
\end{proof}

\medskip

By Theorem~\ref{bigdim} we may assume that $\dim W_j< n-i-1$. Then by Theorem~\ref{general} and its proof (with $\{y_1,\dots,y_m\}=\{x_0,\dots,x_{n-i-1}\}$), there exist $t$, with $0<t\leq \dim W_j$, linear forms $L_1,\dots,L_t$,
which are part of a basis of $W_j$, and linearly independent linear forms $h_1,\dots,h_{n-i-1-t}$,
in variables $x_0,\dots,x_{n-i-1}$, but not involved in $\{L_1,\dots,L_t\}$, such that
\begin{equation}\label{conclusion2}(L_1,\dots,L_t)(h_1,\dots,h_{n-i-1-t})\subset I.\end{equation}
By construction we may assume that $\{L_1,\dots,L_t\}\subset {\mathcal L}^j$. If $1\leq p\leq n-i-1-t$, then $h_p$ ``contributes" the variable $x_{l_p}$; that is, the coefficient
of $x_{l_p}$ in $h_p$ is not zero, and $l_p\neq l_q$ if $p\neq q$.

Recall that $x_j(L_1,\dots,L_t)\subset I$.
If $$(L_1,\dots,L_t)(h_1,\dots,h_{n-i-1-t},x_{n-i},\dots,x_n)\subset I,$$ then the conclusion of Theorem~\ref{main} holds.

As in the proof of Theorem~\ref{bigdim}, let $A_0=\{x_{n-i},\dots,x_n\}$, and let $$V_1=\{x_u\in A_0 |\ x_uL^j_{ab}\notin I\ {\rm for\  some}\ L^j_{ab}\in \{L_1,\dots,L_t\} \}.$$
We may assume that $V_1\neq \emptyset$.
Since $F_{abu}=x_uL^u_{ab}$, by Corollary~\ref{xjGxfG} we have that
$x_u\in V_1$ yields a linear form $H_{c_1u}$, $c_1\in \{a,b\}$, with coefficient of $x_u$ different from zero, such that $F_{c_1uj}=x_jH_{c_1u}$.
Let $H_{V_1}$ be the set of such linear forms $H_{c_1u}$.

Fix $u\in V_1$. By Lemma~\ref{relquad} (2) we have that for all $s=0,\dots,n-i-1$, $$F_{suj}=x_j((-1)^{c_1+s}\mu_{c_1uj}x_s+\nu_{suj}x_u).$$

If $\mu_{c_1uj}\neq 0$, then $H_{l_1u},\dots,H_{l_{n-i-1-t}u}$ are linearly independent, and $$x_j(L_1,\dots,L_t, H_{V_1}, H_{l_1u},\dots,H_{l_{n-i-1-t}u})\subset I.$$
Let $H_{V_1'}=\{H_{V_1}, H_{l_1u},\dots,H_{l_{n-i-1-t}u}\}$, and
let $A_1=A_0\setminus V_1$.
If $(A_1)(H_{V_1'})\subset I$, the conclusion of Theorem~\ref{main} follows.

%Let $x_u\in V_1$ and $x_v\in A_1$. By Lemma~\ref{relquad} (2) we have $F_{c_1uv}\in (x_v)$ if and only if $F_{l_puv}\in (x_v)$ for some $1\leq p\leq n-i-1-t$.
Notice that if $\{u,v\}\subset \{n-i,\dots,n\}$, then by Lemma~\ref{relquad} (2) we have $F_{c_1uv}\in (x_v)$ if and only if $F_{l_puv}\in (x_v)$ for some $1\leq p\leq n-i-1-t$.
Now we proceed as in the proof of Theorem~\ref{bigdim}, and the conclusion of Theorem~\ref{main} holds.

\medskip

We may assume that $\mu_{c_1uj}=0$ for all $x_u\in V_1$. Then $F_{c_1uj}=\nu_{c_1uj}x_jx_u$, so that $x_jx_u\in I$.
If $$(A_1,h_1,\dots,h_{n-i-1-t})(L_1,\dots,L_t,V_1)\subset I,$$ then the conclusion of Theorem~\ref{main} follows.

\medskip

Section~\ref{monomialtimesh} below shows that $(V_1)(h_1,\dots,h_{n-i-1-t})\subset I$. Therefore we may assume that $x_ux_v\notin I$ for some $x_u\in V_1$ and some $x_v\in A_1$; that is, the set $V_2=\{x_v\in A_1 |\ x_vH_{c_1u}\notin I\ {\rm for\  some}\ H_{c_1u}\in H_{V_1} \}$ is not empty.
%By Corollary~\ref{xjGxfG} we have that $F_{c_1vj}\in (x_j)$ and $F_{uvj}\in (x_j)$.
As in the proof of Theorem~\ref{bigdim}, we have that $x_v\in V_2$ yields a linear form $H_{c_2v}$, $c_2\in \{c_1,u\}$, with coefficient of $x_v$ different from zero, such that $F_{c_2vj}=x_jH_{c_2v}$.
Let $H_{V_2}$ be the set of such linear forms.

Fix $v\in V_2$. If $c_2=c_1\in \{a,b\}$, Lemma~\ref{relquad} (2) implies that for all $s=0,\dots,n-i-1$, $F_{svj}=x_j(\pm \mu_{c_2vj}x_s+\nu_{svj}x_v)$. If $\mu_{c_2vj}\neq 0$, then $H_{l_1v},\dots,H_{l_{n-i-1-t}v}$ are linearly independent, and $$x_j(L_1,\dots,L_t, H_{V_1}, H_{V_2}, H_{l_1v},\dots,H_{l_{n-i-1-t}v})\subset I.$$
Let $H_{V_2'}=\{H_{V_2}, H_{l_1v},\dots,H_{l_{n-i-1-t}v}\}$, and let $A_2=A_1\setminus V_2$. If $(A_2)(H_{V_2'})\subset I$, the conclusion of Theorem~\ref{main} follows. Otherwise we proceed as in the proof of Theorem~\ref{bigdim}.

Therefore if $c_2=c_1$ we may assume that $\mu_{c_2vj}=0$. Then $F_{c_2vj}=\nu_{c_2vj}x_jx_v$, so that $x_jx_v\in I$. If $c_2=u$, then $F_{uvj}=x_j(\mu_{uvj}x_u+\nu_{uvj}x_v)$, with $\nu_{uvj}\neq 0$. Since $x_jx_u\in I$, we have that $x_jx_v\in I$.  Furthermore, $F_{c_1vj}=x_j(\mu_{c_1vj}x_{c_1})$. As above, if $\mu_{c_1vj}\neq 0$, then $H_{l_1v},\dots,H_{l_{n-i-1-t}v}$ are linearly independent, and we proceed as in the proof of Theorem~\ref{bigdim}.

Hence we may assume that the condition ``$x_ux_v\notin I$ for some $x_u\in V_1$ and some $x_v\in A_1$'' yields $x_jx_v\in I$. We say that $x_v$ is {\it introduced from} $x_u$.

Section~\ref{monomialtimesh} shows that $(V_2)(h_1,\dots,h_{n-i-1-t})\subset I$. If $$(A_2,h_1,\dots,h_{n-i-1-t})(L_1,\dots,L_t,V_1,V_2)\subset I,$$ the conclusion of Theorem~\ref{main} follows. Otherwise we repeat the argument. Proceeding in this way, at step $l\geq 1$ either we conclude as in the proof of Theorem~\ref{bigdim}, or we introduce a new set of monomials $V_{l}$ such that $x_j(V_l)\subset I$ and $(V_l)(h_1,\dots,h_{n-i-1-t})\subset I$.
Furthermore, by inductively applying (\ref{defg1}), we have that if $x_z\in V_{l}$, then $F_{c_1zj}\in (x_j)$. Therefore we assume that $\mu_{c_1zj}=0$; that is, in the notation of Lemma~\ref{relquad} (2), $P_{zj}=0$.
%otherwise we conclude as in Theorem~\ref{bigdim}.
%because $\lambda_{c_1zj}=\mu_{c_1zj}=0$.
This procedure has to terminate in a finite number of steps, and so the conclusion of Theorem~\ref{main} holds.

\subsection{}\label{monomialtimesh}
To conclude the proof of Theorem~\ref{main} we need to show that for each $l\geq 1$, the set $V_l$ has the property that $(V_l)(h_1,\dots,h_{n-i-1-t})\subset I$.

\medskip

We follow the proof of Theorem~\ref{general} and we consider the explicit description of $h_1,\dots,h_{n-i-1-t}$ in (\ref{conclusion2}).
Recall that $\{y_1,\dots,y_m\}=\{x_0,\dots,x_{n-i-1}\}$.

If $0\leq a<b\leq n-i-1$ and $n-i \leq u \leq n$, we have that $F_{abu}=x_uL^u_{ab}$. Let $$L^u_{ab}=\lambda_{ab}^ux_a+\mu_{ab}^ux_b.$$

\medskip

First we summarize some general facts that will be used often.

 \begin{remark}\label{often} {\rm Let $0\leq a<b\leq n-i-1$, and $n-i\leq d<e\leq n$. By Lemma~\ref{relquad} (3) we have that
$\lambda^d_{ab}=(-1)^{a-e}\nu_{bde}+(-1)^{d-e}\lambda^e_{ab}$ and $\mu^d_{ab}=(-1)^{b-1-e}\nu_{ade}+(-1)^{d-e}\mu^e_{ab}$.}\end{remark}

\begin{lemma}\label{mongeneral} Let $\{c,f,g\}\subset \{0,\dots,n-i-1\}$. Let $T_{fg}$ be a linear form in $x_f$ and $x_g$. Suppose that $x_cT_{fg}\in I$, and that $L^j_{cf}, L^j_{cg}$ are monomials in $x_c$. Let $u\in \{n-i,\dots,n\}\setminus\{j\}$ and suppose that the coefficient of $x_ux_j$ in $F_{cuj}$ is not zero. Then $x_uT_{fg}\in I$.\end{lemma}
%$\nu_{auj}\neq 0$.
 %$x_uL^j_{ab}\notin I$.
\begin{proof} Assume $c<f<g<u<j$, so that $\nu_{cuj}\neq 0$. By Remark~\ref{often} we have that $\mu^u_{cf}=\pm \nu_{cuj}$, $\mu^u_{cg}=\pm \nu_{cuj}$. If $x_uT_{fg}\notin I$, there exists a point $E$ such that $x_u(E)\neq 0$, $T_{fg}(E)\neq 0$, and $x_c(E)= 0$. Without loss of generality assume that $x_f(E)\neq 0$. Recall that  $F_{cfu}=x_u(\lambda^u_{cf}x_c+\mu^u_{cf}x_f)$. Then $F_{cfu}(E)=0$ implies that $\mu^u_{cf}=0$, a contradiction.
 \end{proof}

\begin{lemma}\label{connectedmonomials} Let $u\in \{n-i,\dots,n\}$ and assume that $x_uL^j_{ab}\notin I$ for some $L^j_{ab}\in \{L_1,\dots,L_t\}$. Let $T_{fg}$ be a linear form in $x_f$ and $x_g$, where $\{f,g\}\subset \{0,\dots,n-i-1\}$. Assume that $x_aT_{fg}\in I$ and $x_bT_{fg}\in I$, that $L^j_{af}$, $L^j_{ag}$ are monomials in $x_a$, and $L^j_{bf}$, $L^j_{bg}$ are monomials in $x_b$. Suppose $x_z\in V_l$ is introduced inductively from $x_u\in V_1$. Then $x_zT_{fg}\in I$.\end{lemma}
\begin{proof} We proceed by induction on $l$. If $l=1$, then $z=u$. Since $x_jL^j_{ab}\in I$ and $x_uL^j_{ab}\notin I$, by Corollary~\ref{xjGxfG} we have that $(\nu_{auj},\nu_{buj})\neq (0,0)$. Then $x_uT_{fg}\in I$ by Lemma~\ref{mongeneral}. Now suppose that $l>1$ and that $x_z\in V_l$ is introduced because $x_wx_z\notin I$ for some $x_w\in V_{l-1}$; that is, $x_zH_{c_{l-1}w}\notin I$ for some $H_{c_{l-1}w}\in H_{V_{l-1}}$. Assuming $c_{l-1}<w<z<j$, by construction we have that $(\nu_{c_{l-1}zj},\nu_{wzj})\neq (0,0)$. If $c_{l-1}=c_1\in \{a,b\}$ and $\nu_{c_1zj}\neq 0$, we conclude by Lemma~\ref{mongeneral}. Otherwise $\nu_{c_lzj}\neq 0$, where $x_{c_l}\in V_p$ for some $p<l$, so that $x_{c_l}T_{fg}\in I$. Assume by contradiction that $x_zT_{fg}\notin I$. By Lemma~\ref{xgGxfG} and Lemma~\ref{relquad} (2) we have that $\mu_{fc_lz}=\mu_{gc_lz}=\mu_{c_1c_lz}=0$. Recall that we are assuming $\mu_{c_1c_lj}=0$. Then by (\ref{defg3}) we have that $\nu_{c_lzj}=0$, a contradiction.
\end{proof}

\medskip

Now suppose that (\ref{conclusion2}) is given by (\ref{eqbadblock}) of Remark~\ref{badblock},
$$(B^{12})(Y_C)\subset I.$$ Let $x_u\in V_1$ be such that $x_uL^j_{ab}\notin I$. By construction we have that $x_a,x_b\in B^{12}$. If $s\in C$ let $T_s= x_s$. Then by Lemma~\ref{connectedmonomials}, we have that $(V_l)(Y_C)\subset I$, as desired.

\medskip

Next we suppose that (\ref{conclusion2}) is given by (\ref{Gak}) of Lemma~\ref{monomial},
$$(B^{12}, Y_{M_k}, Y_{C_k})(G_{A_k})\subset I,$$ for some $k\geq 1$. For simplicity of notation we may assume that $x_1=y_1$ and $x_2=y_2$. Recall that $G_{A_k}$ consists of forms $G_{fg}$ defined in (\ref{defG}), $$G_{fg}=((-1)^{f-1}\lambda^j_{1g}-\lambda^j_{fg})x_f+((-1)^g\lambda^j_{1f}-\mu^j_{fg})x_g.$$ Therefore (\ref{Gak}) includes equations of type (\ref{BGN}).

Let $x_u\in V_1$. By Lemma~\ref{connectedmonomials}, it suffices to consider the cases $x_uL^j_{12}\notin I$, where $L^j_{12}=\lambda^j_{12}x_1$, and $x_uL^j_{2q}\notin I$, with $x_q\in Y_{M_k}$; that is, $L^j_{2q}=(-1)^q \lambda^j_{12}x_q$.

\medskip

Suppose $x_uL^j_{12}\notin I$. Let $z\in V_l$, $l\geq 1$, be obtained inductively from $x_u$. We have that $F_{fgz}=x_z(\lambda^z_{fg}x_f+\mu^z_{fg}x_g).$
By Remark~\ref{often}, assuming that $f<g<z<j$, we have that $$\lambda^z_{fg}=(-1)^{f-j}\nu_{gzj}+(-1)^{z-j}\lambda^j_{fg}$$
and $$\mu^z_{fg}=(-1)^{g-1-j}\nu_{fzj}+(-1)^{z-j}\mu^j_{fg}.$$ Since $x_1G_{fg}\in I$, if $\nu_{1zj}\neq 0$, we have that $x_zG_{fg}\in I$ by Lemma~\ref{mongeneral}. Therefore we may assume that $\nu_{1zj}=0$.

%First we show that $x_uG_{fg}\in I$. We have that $F_{fgu}=x_u(\lambda^u_{fg}x_f+\mu^u_{fg}x_g).$
%By Remark~\ref{often}, assuming that $1<f<g<u<j$ we have that $\lambda^u_{fg}=(-1)^{f-j}\nu_{guj}+(-1)^{u-j}\lambda^j_{fg}$
%and $\mu^u_{fg}=(-1)^{g-1-j}\nu_{fuj}+(-1)^{u-j}\mu^j_{fg}$.

First we show that $x_uG_{fg}\in I$.
%Since $x_1G_{fg}\in I$, if $\nu_{1uj}\neq 0$, we have that $x_uG_{fg}\in I$ by Lemma~\ref{mongeneral}.
Since $\nu_{1uj}= 0$, by Remark~\ref{often} we have that $\mu^u_{1f}=\mu^u_{1g}=0$. Then $F_{1fu}=F_{1gu}=0$, since $x_1x_u\notin I$. Then by Remark~\ref{often} we have that $\nu_{fuj}=(-1)^{u}\lambda^j_{1f}$ and $\nu_{guj}=(-1)^{u}\lambda^j_{1g}$.
It follows that $$F_{fgu}=(-1)^{u-j+1}x_uG_{fg},$$ and so $x_uG_{fg}\in I$.

\medskip

Now suppose $x_v\in V_2$ is introduced because $x_ux_v\notin I$.
%If $\nu_{1vj}\neq 0$, we have that $x_vG_{fg}\in I$ by Lemma~\ref{mongeneral}.
If $\mu_{1uv}\neq 0$, we have that $x_vG_{fg}\in I$ by Lemma~\ref{xgGxfG} and Lemma~\ref{relquad} (2). Therefore we may assume that $\mu_{1uv}=0$.
Since $x_1x_v\in I$ and $x_1x_u\notin I$, we have that $\lambda_{1uv}=0$. Hence for all $s=0,\dots,n-i-1$, we have that $F_{suv}=\nu_{suv}x_ux_v$. It follows that $\nu_{suv}=0$. In particular $\nu_{fuv}=\nu_{guv}=\nu_{1uv}=0$. Then by (\ref{defg4}) we have that $\nu_{fvj}=(-1)^{v-u}\nu_{fuj}$, $\nu_{gvj}=(-1)^{v-u}\nu_{guj}$, and $\nu_{1uj}=\nu_{1vj}=0$. As above, $\nu_{1uj}=0$ implies that $\nu_{fuj}=(-1)^{u}\lambda^j_{1f}$ and $\nu_{guj}=(-1)^{u}\lambda^j_{1g}$.
%Furthermore, by Remark~\ref{often} we have that $\lambda^v_{fg}=(-1)^{f-j}\nu_{gvj}+(-1)^{v-j}\lambda^j_{fg}$
%and $\mu^v_{fg}=(-1)^{g-1-j}\nu_{fvj}+(-1)^{v-j}\mu^j_{fg}$.
It follows that $F_{fgv}=(-1)^{v-j+1} x_vG_{fg},$ and so $x_vG_{fg}\in I$.

\medskip

We proceed by induction on $l$. Suppose that $l\geq 3$, and that $x_z=x_{u_l}$ is introduced inductively from $x_u=x_{u_1}$. Since $x_{u_p}G_{fg}\in I$ for all $1\leq p<l$, by Lemma~\ref{xgGxfG} we may assume that $\mu_{1u_pu_l}= 0$ for all $1\leq p<l$.

Let $d<e<f$. Observe that if $x_dx_f\in I$ and $x_ex_f\notin I$, then $\nu_{def}=0$. If $x_dx_f\in I$ and $x_dx_e\notin I$, then $\lambda_{def}=0$. It follows that $\lambda_{1u_1u_p}=0$ for all $1<p\leq l$, and $\lambda_{u_{p-1}u_pu_m}=0$ for all $1<p< l$ and $p<m\leq l$. By inductively applying (\ref{defg1}), we have that $\lambda_{1u_{p-1}u_p}=0$ for all $1<p\leq l$.

Similarly, applying (\ref{defg3}), it follows that $\mu_{1u_{p-1}u_p}=0$ for all $1<p\leq l$.
%Inductively applying this observation and Remark~\ref{defg}, it follows that $\lambda_{1u_{p-1}u_p}=\mu_{1u_{p-1}u_p}= 0$ for all $1<p\leq l$.
Since $x_{u_{p-1}}x_{u_p}\notin I$, we have that $\nu_{su_{p-1}u_p}=0$ for all $1<p\leq l$ and $0\leq s\leq n-i-1$. Then by (\ref{defg4}) we have that for all $1<p\leq l$ and $0\leq s\leq n-i-1$, $\nu_{su_pj}=(-1)^{u_p-u_{p-1}}\nu_{su_{p-1}j}$. It follows that $\nu_{fzj}=(-1)^{z-u}\nu_{fuj}$, $\nu_{gzj}=(-1)^{z-u}\nu_{guj}$, and $\nu_{1uj}=\nu_{1zj}=0$.
Hence $\nu_{fuj}=(-1)^{u}\lambda^j_{1f}$ and $\nu_{guj}=(-1)^{u}\lambda^j_{1g}$.
%Furthermore, by Remark~\ref{often} we have that $\lambda^v_{fg}=(-1)^{f-j}\nu_{gvj}+(-1)^{v-j}\lambda^j_{fg}$
%and $\mu^v_{fg}=(-1)^{g-1-j}\nu_{fvj}+(-1)^{v-j}\mu^j_{fg}$.
It follows that $F_{fgz}=(-1)^{z-j+1} x_zG_{fg},$ and so $x_zG_{fg}\in I$.

\medskip

Now suppose that $x_uL^j_{2q}\notin I$. As in the proof of Lemma~\ref{monomial} we assume that $f<g<q<u<j$. Recall that $q$ has the property that the coefficient of $x_q$ in $G_{fq}$ and in $G_{gq}$ is zero; that is $\mu^j_{fq}=(-1)^q\lambda^j_{1f}$, and $\mu^j_{gq}=(-1)^q\lambda^j_{1g}$.

Let $x_u\in V_1$. As before we may assume that $\nu_{quj}=0$, and so $\lambda^u_{fq}=\lambda^u_{gq}=0$. Then $F_{fqu}=F_{gqu}=0$, since $x_qx_u\notin I$. It follows that $\nu_{fuj}=(-1)^{u-q}\mu^j_{fq}=(-1)^{u}\lambda^j_{1f}$ and $\nu_{guj}=(-1)^{u}\lambda^j_{1g}$.
Then $F_{fgu}=(-1)^{u-j+1}x_uG_{fg},$ and so $x_uG_{fg}\in I$.

If $l\geq 2$ we repeat the proof of the previous case, with $x_q$ instead of $x_1$, and we obtain that $x_zG_{fg}\in I$.
%As before, $x_pG_{fg}\in I$, $x_uG_{fg}\notin I$, and $x_px_u\notin I$ yields $F_{fpu}=F_{gpu}=0$. Then $\nu_{fuj}=(-1)^{u-p}\mu^j_{fp}$ and $\nu_{guj}=(-1)^{u-p}\mu^j_{gp}$. Since $\mu^j_{fp}=(-1)^p\lambda^j_{1f}$, and $\mu^j_{gp}=(-1)^p\lambda^j_{1g}$, as above we conclude that $x_uG_{fg}\in I$, a contradiction.

\medskip

Now we only need to consider the case when (\ref{conclusion2}) is given by (\ref{final}), $$(B_1,B_{d_2},\dots,B_{d_r},D_{j_1})(V_N,V_L\setminus V_{D_{j_1}},B_{d_{r+1}},\dots,B_{d_k})\subset I.$$
Notice that (\ref{final}) includes the cases (\ref{basis+varout}), (\ref{eqonlybin}), and (\ref{onemonblock}).
%when $V_N=\emptyset$, $k=0$ (that is, the basis consists only of binomials), $k=1$, and $r=1$.
Here the linear forms $h_1,\dots,h_{n-i-1-t}$ are all monomials.

We will need the following observations.

\begin{remark}\label{xaxpxbxpinI} {\rm Let $u\in \{n-i,\dots,n\}$ be such that $x_uL^j_{ab}\notin I$ for some $L^j_{ab}\in \{L_1,\dots,L_t\}$. Let $x_p\in \{h_1,\dots,h_{n-i-1-t}\}$, and assume that $x_ax_p\in I$, $x_bx_p\in I$, that $L^j_{ap}$ is a monomial in $x_a$, and $L^j_{bp}$ is a monomial in $x_b$. Suppose that $x_z\in V_l$, $l\geq 1$, is introduced inductively from $x_u\in V_1$. Then $x_zx_p\in I$, by Lemma~\ref{connectedmonomials}.}
\end{remark}

\begin{lemma}\label{LapLbpzero} Let $u\in \{n-i,\dots,n\}$ be such that $x_uL^j_{ab}\notin I$ for some $L^j_{ab}\in \{L_1,\dots,L_t\}$. Let $x_p\in \{h_1,\dots,h_{n-i-1-t}\}$, and assume that $L^j_{ap}=L^j_{bp}=0$. Suppose that $x_z\in V_l$, $l\geq 1$, is introduced inductively from $x_u\in V_1$. Then $F_{pzj}\in (x_p)$, and $x_zx_p\in I$.
\end{lemma}
\begin{proof} Assume that $a<b<p<z<j$. Let $x_z\in V_l$, $l\geq 1$. By Remark~\ref{often} we have that $\lambda^z_{ap}=\pm\nu_{pzj}$, $\lambda^z_{bp}=\pm\nu_{pzj}$, $\mu^z_{ap}=\pm\nu_{azj}$, and $\mu^z_{bp}=\pm\nu_{bzj}$. We show by induction on $l$ that $\nu_{pzj}=0$ and that $x_px_z\in I$.

If $l=1$, then $z=u$, and $\nu_{puj}=0$ by Lemma~\ref{mongeneral}. It follows that $\lambda^u_{ap}=\lambda^u_{bp}=0$, and so $F_{apu}=\mu^u_{ap}x_ux_p$, and $F_{bpu}=\mu^u_{bp}x_ux_p$. Since $x_jL^j_{ab}\in I$ and $x_uL^j_{ab}\notin I$, by Corollary~\ref{xjGxfG} we have that $(\nu_{auj},\nu_{buj})\neq (0,0)$. Therefore $(\mu^u_{ap}, \mu^u_{bp})\neq (0,0)$ and $x_px_u\in I$.

Now suppose that $l>1$ and that $x_z\in V_l$ is introduced because $x_wx_z\notin I$ for some $x_w\in V_{l-1}$; that is, $x_zH_{c_{l-1}w}\notin I$ for some $H_{c_{l-1}w}\in H_{V_{l-1}}$. Assuming $c_{l-1}<w<z<j$, by construction we have that $(\nu_{c_{l-1}zj},\nu_{wzj})\neq (0,0)$.

Now $x_px_w\in I$, and $x_wx_z\notin I$ imply that $\nu_{pwz}=0$. Since $\nu_{pwj}=0$ by the induction hypothesis, it follows from (\ref{defg4}) that $\nu_{pzj}=0$.
Then $\lambda^z_{ap}=\lambda^z_{bp}=0$, and so $F_{apz}=\mu^z_{ap}x_zx_p$, and $F_{bpz}=\mu^z_{bp}x_zx_p$.

If $c_{l-1}=c_1\in \{a,b\}$ and $\nu_{c_1zj}\neq 0$, then $(\mu^z_{ap}, \mu^z_{bp})\neq (0,0)$ and $x_zx_p\in I$.
Otherwise $\nu_{c_lzj}\neq 0$, where $x_{c_l}\in V_p$ for some $p<l$, so that $x_{c_l}x_p\in I$. If $x_zx_p\notin I$, then we have that $\mu_{pc_lz}=\mu_{c_1c_lz}=0$. Recall that we are assuming $\mu_{c_1c_lj}=0$. Then by (\ref{defg3}) we have that $\nu_{c_lzj}=0$, a contradiction.
\end{proof}

\begin{lemma}\label{connected2} Let $u\in \{n-i,\dots,n\}$ and $\{p,q\}\subset \{0,\dots,n-i-1\}$. Suppose that $x_ux_p\in I$, that $F_{puj}\in (x_p)$, and that the coefficient of $x_q$ in $L^j_{pq}$ is not zero. Then $x_ux_q\in I$.
\end{lemma}
\begin{proof} Assume $p<q<u<j$. By Remark~\ref{often} we have that $\mu^u_{pq}=\pm \mu^j_{pq}\neq 0$. Then $F_{pqu}=x_u(\lambda^u_{pq}x_p+\mu^u_{pq}x_q)$ and $x_ux_p\in I$ imply that $x_ux_q\in I$.
\end{proof}

\begin{lemma}\label{connected3} Let $u\in \{n-i,\dots,n\}$ be such that $x_uL^j_{ab}\notin I$ for some $L^j_{ab}\in \{L_1,\dots,L_t\}$. Let $x_p\in \{h_1,\dots,h_{n-i-1-t}\}$, and assume that $L^j_{ap}$ and $L^j_{bp}$ are monomials in $x_p$. Suppose that $x_zx_p\in I$ for all $x_z$ introduced inductively from $x_u\in V_1$. If the coefficient of $x_q$ in $L^j_{pq}$ is not zero, then $x_zx_q\in I$.
\end{lemma}
\begin{proof} By Lemma~\ref{mongeneral} we have that $\nu_{puj}=0$. Now the proof of Lemma~\ref{LapLbpzero} shows that $\nu_{pzj}=0$.
%$\nu_{puj}=0$ and $x_zx_p\in I$ for all $x_z$ introduced from $x_u$ implies that $\nu_{pzj}=0$.
Then by Lemma~\ref{connected2} we have that $x_zx_q\in I$.
\end{proof}

\medskip

We are now ready to conclude the proof that $(V_l)(h_1,\dots,h_{n-i-1-t})\subset I$ for all $l\geq 1$. Let $u\in \{n-i,\dots,n\}$ be such that $x_uL^j_{ab}\notin I$ for some $L^j_{ab}\in \{L_1,\dots,L_t\}$. Suppose that $x_z\in V_l$, $l\geq 1$, is introduced inductively from $x_u\in V_1$. Let $x_p\in \{V_N, V_L\setminus V_{D_{j_1}},B_{d_{r+1}},\dots,B_{d_k}\}$. Let $x_{a_{d_{r+1}}}, \dots, x_{a_{d_k}}$ be generators of $B_{d_{r+1}},\dots,B_{d_k}$ respectively.

\medskip

First assume that $L^j_{ab}\in D_{j_1}$. If $x_p \in V_N \cup (V_L\setminus V_{D_{j_1}})\cup \{x_{a_{d_{r+1}}}\}\cup \dots \cup \{x_{a_{d_k}}\}$, then $x_zx_p\in I$ by Lemma~\ref{LapLbpzero}, since $L^j_{ap}=L^j_{bp}=0$.
If $x_p\in \{B_{d_{r+1}},\dots,B_{d_k}\}$ is not a generator of one of the blocks, then $x_zx_p\in I$ by inductively applying Lemma~\ref{connected3}.

\medskip

Next assume that $L^j_{ab}$ is one of the generators of $B_1,B_{d_2},\dots,B_{d_r}$. Recall that $x_b\in V_{D_{j_1}}$. If $x_p\in V_N \cup (V_L\setminus V_{D_{j_1}})$, we may assume that $L^j_{ap}=0$, otherwise by Lemma~\ref{monomial} and by the proof of Lemma~\ref{monandbin}, we can reduce to the previous case given by equation (\ref{Gak}). If $x_p\in \{x_{a_{d_{r+1}}}\}\cup \dots \cup \{x_{a_{d_k}}\}$ we have that $L^j_{ap}=0$ by Remark~\ref{2blocks}. We also have that $L^j_{bp}=0$ if $x_p \in V_N \cup (V_L\setminus V_{D_{j_1}})\cup \{x_{a_{d_{r+1}}}\}\cup \dots \cup \{x_{a_{d_k}}\}$. Then by Lemma~\ref{LapLbpzero} we have that $x_zx_p\in I$.

If $x_p\in \{B_{d_{r+1}},\dots,B_{d_k}\}$ is not a generator, then $L^j_{ap}$ is a monomial in $x_p$ (otherwise the block containing $x_a$ would not be maximal), and $L^j_{bp}$ is a monomial in $x_p$.  Then $x_zx_p\in I$ by inductively applying Lemma~\ref{connected3}.

\medskip

Last assume that $x_a$ and $x_b$ belong to one of the blocks $B_1,B_{d_2},\dots,B_{d_r}$. If $x_p \in V_N \cup (V_L\setminus V_{D_{j_1}})\cup \{x_{a_{d_{r+1}}}\}\cup \dots \cup \{x_{a_{d_k}}\}$, then $L^j_{ap}$ and $L^j_{bp}$ are monomials in $x_a$ and $x_b$ respectively, and so by Remark~\ref{xaxpxbxpinI} we have that $x_zx_p\in I$. The same holds if $x_p$ is not a generator of $B_{d_{r+1}},\dots,B_{d_k}$, but the block containing $x_p$ has been constructed after the block containing $x_a$ and $x_b$.

If the block containing $x_p$ has been constructed before the block containing $x_a$ and $x_b$, then $L^j_{ap}$ and $L^j_{bp}$ are monomials in $x_p$ and by inductively applying Lemma~\ref{connected3} we have that $x_zx_p\in I$.

\medskip

This concludes the proof of $(V_l)(h_1,\dots,h_{n-i-1-t})\subset I$, of Theorem~\ref{main}, and of Conjecture~\ref{theconjecture}.

\end{proof}

\end{document}